\documentclass[11pt]{amsart}
\usepackage[latin1]{inputenc}
\usepackage[english]{babel}
\usepackage{amsmath,amssymb,amsthm,amsfonts}
\usepackage[all]{xy}
\usepackage{tikz}
\usepackage{tikz-cd}

\usepackage{latexsym,dsfont}
\usepackage{enumerate}
\usepackage{graphicx}
\usepackage{latexsym}

\makeatletter
\@namedef{subjclassname@2010}{ \textup{2010} Mathematics Subject Classification}
\makeatother

\newtheorem{theorem}{Theorem}

\newtheorem{lemma}{Lemma}
\newtheorem{corollary}[lemma]{Corollary}
\newtheorem{proposition}[lemma]{Proposition}

\theoremstyle{definition}

\newtheorem{remark}[lemma]{Remark}

\newtheorem{problem}[lemma]{Problem}

\numberwithin{equation}{section}

\numberwithin{equation}{section}
\DeclareMathOperator{\supp}{supp}
\DeclareMathOperator{\sign}{sign}

\begin{document}
\title{Monotone substochastic operators and a new Calder\'on couple}
\author[Le\'snik]{Karol Le\'snik}
\address[Karol Le{\'s}nik]{Institute of Mathematics\\
Pozna\'n University of Technology, ul. Piotrowo 3a, 60-965 Pozna{\'n}, Poland}
\email{\texttt{klesnik@vp.pl}}


\begin{abstract}
An important result on submajorization, which goes back to Hardy, Littlewood and P\'olya, states that $b\preceq a$ if and only if there is a doubly stochastic matrix $A$ such that $b=Aa$. We prove that under monotonicity assumptions on vectors $a$ and $b$ respective matrix $A$ may be chosen monotone. 
This result is then applied to show that  $(\widetilde{L^p},L^{\infty})$ is a Calder\'on couple for $1\leq p<\infty $, where $\widetilde{L^{p}}$ is the K\"othe dual of the Ces\`aro space $Ces_{p'}$ (or equivalently the down space $L^{p'}_{\downarrow}$). In particular, $(\widetilde{L^1},L^{\infty})$ is a Calder\'on couple and this complements the result of \cite{MS06} where it was shown that $(L^{\infty}_{\downarrow},L^{1})$ is a Calder\'on couple.      
\end{abstract}
\maketitle

\footnotetext[1]{2010 \textit{Mathematics Subject Classification}: 46E30, 46B42, 46B70}
\footnotetext[2]{\textit{Key words and phrases}: substochastic operator, doubly stochastic matrix, Calder\'on couple}



\section{Introduction}

The classical Hardy, Littlewood and P\'olya theorem on submajorization states that for two vectors $a,b\in \mathbb{R}^n$, $a\preceq b$ is equivalent with existence of a doubly stochastic matrix $A$ satisfying $Ab=a$ (see for example \cite{BS88}, \cite{MOA11} or \cite{Mi88}). These ideas were developed by Calder\'on in \cite{Ca66} (cf. \cite{Ry65}) where he proved that for two functions $f,g\in L^0$ such that $f\prec g$ and $g\in L^1+L^{\infty}$ there is a substochastic operator $T$  such that $Tg=f$. This allowed him to complete the ``only if" part of his famous theorem which states that a space $X$ is interpolation space between $L^1$ and $L^{\infty}$ if and only if $X$ satisfies the following property: if two functions $f,g\in L^0$ are such that $f\prec g$ and $g\in X$, then also $f\in X$ (see \cite{Ca66}). The ``only if'' part seems to have deeper impact and to be more spectacular. Moreover, it initiated investigations of the so called Calder\'on couples and $(L^1,L^{\infty})$ was the first such a couple. Recall that $(X_0,X_1)$ is a {\it Calder\'on couple} if each interpolation space between $X_0$ and $X_1$ is $K$-monotone (with respect to $(X_0,X_1)$).    
Thanks to the $K$-divisibility theorem of Brudnyi-Krugljak each $K$-monotone space may be represented by the $K$-method of interpolation, which means that all interpolaion spaces for a Calder\'on couple may be produced by the $K$-method. Some of the most important papers in interpolation theory deal with the problem of being a Calder\'on couple, see for example \cite{AC84}, \cite{Cw76}, \cite{Ka93}, \cite{LS71} and \cite{Sp78}.

We have sketched shortly how Hardy-Littlewood-P\'olya theorem evolved into Calder\'on couples, because in the paper we will adopt the mentioned steps in such a way that we are able to prove that  $(\widetilde{L^1},L^{\infty})$ is a Calder\'on couple. We start with proving the monotone version of the Hardy-Littlewood-P\'olya theorem and apply it to provide respective, monotone version of Calder\'on's theorem. Notice, that actually the second one (and so also the first one) was already proved by Bennett and Sharpley in \cite{BS86}.  They used it in an alternative prove of the K-divisibility theorem for Gagliardo couples. Our method is however essentially different and we present it in the first part of the paper. 

The second part is devoted to proving that $(\widetilde{L^1},L^{\infty})$ is a  Calder\'on couple, where $\widetilde{L^p}$ is defined as 
\[
\widetilde{L^p}=\{f\in L^0:\widetilde{f}\in L^p\}
\]
for $\widetilde{f}$ being the nonincreasing majorant of $f$, i.e. 
\begin{equation}
\widetilde{f}(t)={\rm ess}\sup_{s\geq t} |f(s)|,\ t>0. \label{deffalki}
\end{equation}
In \cite{MS06} it was shown that  $(L^{\infty}_{\downarrow},L^{1})$ is a Calder\'on couple, where $L^{\infty}_{\downarrow}$ means a down space of $L^{\infty}$. On the other hand, the K\"othe dual of a down space  $L^{p}_{\downarrow}$ is $\widetilde{L^{p'}}$, so the couple $(\widetilde{L^1},L^{\infty})$ is just the dual one to $(L^{\infty}_{\downarrow},L^{1})$. In the remaining part we discuss a possible extensions of this result. In particular, we adopt Dimitriev's results to show that $(\widetilde{\Lambda_{\varphi}},L^{\infty})$ is a relative Calder\'on couple with respect to $(\widetilde{L^1},L^{\infty})$. Moreover, using Avni-Cwikel theorem \cite{AC12} we conclude that $(\widetilde{L^p},L^{\infty})$ is a Calder\'on couple also for $1<p<\infty$.

\section{Basic definitions}
Let $\mu$ be a Lebesgue measure on $\mathbb{R}_+$ and denote by $L^0$ the space of all (equivalence classes of) real-valued functions on $\mathbb{R}_+$. 
A Banach space $X\subset L^{0}$ is called a {\it Banach function space} if $f\in X,g\in L^{0},|g|\leq |f|$ implies $g\in X$ and $||g||_X\leq
||f||_X$. We will also understand that there is $f\in X$ with $f(t)>0$ for
each $t>0$. 

By a {\it symmetric space} we mean a Banach function space $X$ with the additional property that for any two equimeasurable functions 
$f \sim g, f, g \in L^{0}$ (that is, they have the same distribution functions $d_{f}\equiv d_{g}$, where 
$d_{f}(\lambda) = \mu(\{t>0: |f(t)|>\lambda \}),\lambda \geq 0,$ and $f\in X$ we have $g\in X$ and $\| f\|_{X} = \| g\|_{X}$. In particular, 
$\| f\|_{X}=\| f^{\ast }\|_{X}$, where $f^{\ast }(t)=\mathrm{\inf } \{\lambda >0\colon \ d_{f}(\lambda ) < t\},\ t\geq 0$ is the nonincreasing rearrangement of $f$.
For more informations on Banach function spaces and symmetric spaces we refer to \cite{KPS82} or \cite{BS88}.



For $a=(a_1,...,a_n)\in \mathbb{R}^n$, $a^*$ is the vector produced by permuting entries of $|a|$ in nonincreasing order. Writing $b\prec a$, for $a,b\in \mathbb{R}^n$ we understand that
\[
\sum_{i=1}^kb^*_i\leq \sum_{i=1}^ka^*_i\ {\rm for~ each }\ 0<k\leq n,
\]
while the relation $b\preceq a$ means that $b\prec a$ and additionally 
$$
\sum_{i=1}^nb^*_i= \sum_{i=1}^na^*_i.
$$ 

Given a matrix $A$ we shall not distinguish it with the respective linear operator defined by $A$ and write just $Ax$ when understand it as $x\cdot A^T$.
A positive matrix $A=(a_{ij})_{i,j=1}^n$  (here positivity means that $0\leq a_{ij}$ for all $i,j$, or equivalently $0\leq Aa$ for each $0\leq a\in \mathbb{R}^n$) is called {\it doubly stochastic} when 
\begin{equation}
\sum_{j=1}^na_{ij}= \sum_{j=1}^na_{ji}=1\ {\rm for~ each }\ 0<i\leq n. \label{matrixDS}
\end{equation}
If all above sums are just less or equal one, we say about {\it substochastic} matrix. Equivalently, a positive matrix $A$ is doubly stochastic (substochastic) if and only if $Aa\preceq a$  ($Aa\prec a$) for each $0\leq a\in \mathbb{R}^n$. A positive matrix $A$ will be called {\it monotone} if it is positive and $Aa$ is nonincreasing for each nonincreasing $0\leq a\in \mathbb{R}^n$.
It is easy to see that a positive square matrix $A=(a_{ij})_{i,j=1}^n$ is monotone if and only if for each $k=1,\dots,n$ 
\begin{equation}
\sum_{j=1}^ka_{i,j}\geq \sum_{j=1}^ka_{i+1,j}\ {\rm for~ any }\ 0<i\leq n-1. \label{matrixMP}
\end{equation}
In fact, sufficiency of (\ref{matrixMP}) is a consequence of the Hardy lemma and  necessity comes by applying definition to vectors $e_1=[1,0,...,0], e_1+e_2=[1,1,0,...,0]$, etc. (see for example \cite[Chapter 2.E]{MOA11}). 

We shall need also continuous versions of the above objects. A linear positive operator (in the sense that $0\leq f$ implies $0\leq Tf$) defined on $L^1+L^{\infty}$, mapping continuously $L^1$ into  $L^1$ and  $L^{\infty}$ into  $L^{\infty}$ with both norms less or equal one is called {\it substochastic}. This is equivalent with $Tf\prec f$ for each $f\in L^1+L^{\infty}$, when $T$ is positive (cf. \cite{KPS82}, p. 84). 
To speak about monotone functions in the setting of Banach function spaces we need to precise the notion of monotonicity to make it insensitive to perturbation on a set of measure zero. Hence we will understand that $0\leq f\in L^0$ is nonincreasing if it is nonincreasing (in a classical sense) on some set $S\subset \mathbb{R}_+$ such that $\mu(\mathbb{R}_+\backslash S)=0$.
It will be useful in the sequel to use also equivalent formulation, which says that  $0\leq f\in L^0$ is nonincreasing if for each $x>0$ 
\begin{equation}
{\rm ess}\sup_{s\geq x} |f(s)|\leq{\rm ess} \inf_{s\leq x} |f(s)| \label{defmon}
\end{equation}
(see \cite[Theorem 2.4]{Si94}).
Let $X,Y$ be two Banach function spaces. We shall say that an operator $T:X\mapsto Y$ is {\it monotone} if it is positive  and for each nonincreasing $0\leq f \in X$, $Tf$ is also nonincreasing.

Let us recall that a {\it Banach limit} is a linear functional $\eta \in (l^{\infty})^*$ with the following properties:
\begin{itemize}
\item[(i)] if $\lim_{n\rightarrow \infty} x_n$ exists, then $\lim_{n\rightarrow \infty} x_n=\eta(x_n)$,
\item[(ii)] $x_n\geq 0$ for each $n$ implies $\eta(x_n)\geq 0$,
\item[(iii)] $\eta(x_n)=\eta(x_{n+1})$.
\end{itemize}
Such a limits arise from applying  the Hahn-Banach theorem to the respective subspace of $l^{\infty}$. Alternatively, one can define $\eta$ to be a limit of $([Cx]_n)$ with respect to a nonprimal ultrafilter, where $[Cx]_n=\frac{1}{n}\sum_{i=1}^n x_i$. Notice that we use the convention of writing $\eta(x_n)$  rather than $\eta((x_n)_{n=1}^{\infty})$ just to simplify further notion. In the case of multiple indexes, we shall write $\eta_{k\uparrow}(x_k^n)$ to emphasize that we mean $\eta((x_k^n)_{k=1}^{\infty})$.  


For a given Banach function space $X$ we define the space $\widetilde{X}$ as 
\[
\widetilde{X}=\{f\in L^0:\widetilde{f}\in X\}
\]
with the norm given by 
\[
\|f\|_{\widetilde{X}}=\|\widetilde{f}\|_{X}.
\]
To ensure that such a space is a Banach function space in the sense of our definition, we will assume that for a Banach function space $X$ there is a nonincreasing $f\in X$ with $f(t)>0$ for each $t>0$. Spaces $\widetilde{X}$ appear in a natural way in different contexts. It seems that first time such spaces in general form (for $X$ being symmetric) were defined by
Sinnamon, who proved that they are duals of down spaces, i.e. $(X_{\downarrow})'=\widetilde{X}$ (see \cite{Si94}, \cite{Si01}, \cite{Si03} and \cite{Si07}). On the other hand, the space $\widetilde{L^1}$ has appeared to be the K\"othe dual of Ces\`aro space $Ces_{\infty}$ already in the early paper \cite{KKL48}. Recently, it was also proved in \cite{LM14a} that they are duals of Ces\`aro type spaces even for not necessarily symmetric $X$, i.e. $(CX)'=\widetilde{(X')}$ (for more informations on such a spaces and their history see \cite{LM14a} and references therein, cf. \cite{AM09}, \cite{KMS07}). On the other hand, a space $\widetilde{X}$ is just the space $X(Q)$ associated with the cone of positive nonincreasing functions $Q$ considered in \cite{CC05} and \cite{CEP99}.

For two couples of Banach function spaces $(X_0,X_1),(Y_0,Y_1)$ and a linear operator $T$ acting from $X_0+X_1$ into $Y_0+Y_1$ we write $T: (X_0,X_1)\rightarrow (Y_0,Y_1)$ when $T: X_0\rightarrow Y_0$ and $T: X_1\rightarrow Y_1$ with 
\[
\|T\|_{(X_0,X_1)\rightarrow (Y_0,Y_1)}=\max\{\|T\|_{X_0\rightarrow Y_0},\|T\|_{X_1\rightarrow Y_1}\}<\infty.
\]
For $f\in X_0+X_1$ the $K$-functional of $f$ with respect to the couple $(X_0,X_1)$ is defined as
\[
 K(t,f;X_0,X_1)=\inf\{\|f_0\|_{X_0}+t\|f_1\|_{X_1}:f=f_0+f_1\}~~ {\rm for} ~t>0.
\]

Having two couples $(X_0,X_1),(Y_0,Y_1)$ of Banach function spaces we say that $(X_0,X_1)$ is a {\it relative Calder\'on couple} with respect to $(Y_0,Y_1)$ if for each $f\in X_0+X_1$, $g\in Y_0+Y_1$ the inequality $K(t,g;Y_0,Y_1)\leq K(t,f;X_0,X_1)$ for all $t>0$ implies that there exists $T:(X_0,X_1)\rightarrow (Y_0,Y_1)$ with $Tf=g$. If $(X_0,X_1)=(Y_0,Y_1)$ we say  simply about a {\it Calder\'on couple}. See \cite{AC84}, \cite{BK91}, \cite{BL76}, \cite{BS88}, \cite{Cw76}, \cite{KPS82} and \cite{Sp78} for more informations on interpolation spaces and Calder\'on couples.

\section{Monotone version of Calder\'on and Hardy, Littlewood, P\'olya theorems}

The classical theorem of Calder\'on states that if $g\prec f$, then there is a substochastic operator such that $Tf=g$. It may be regarded as a continuous version of the Hardy, Littlewood and P\'olya  theorem which ensures existence of a doubly stochastic matrix $A$ such that $Aa=b$ provided $b\preceq a$ where $a,b\in \mathbb{R}^n$ (\cite[Theorem 2.7, p. 108]{BS88}, cf. \cite{Mi88} and \cite{MOA11}). 
We should be interested in the following monotone refinements of those theorems.
 
\begin{theorem}[Bennett-Sharpley 1986]\label{CaldMitjMP}
Let $0\leq f,g \in L^1+L^{\infty}$ be both nonincreasing and suppose that $g\prec f$. Then there is a substochastic monotone operator $T$ such that $Tf=g$. 
\end{theorem}

\begin{theorem}\label{thmHLP}
Let $0\leq a,b \in \mathbb{R}^n$ be both nonincreasing. If $b\preceq a$, then there exists a doubly stochastic monotone matrix $A$ such that $Aa=b$. If just $b\prec a$, then the matrix $A$ may be chosen to be substochastic and monotone.   
\end{theorem}

Notice that Theorem \ref{CaldMitjMP} was already proved by Bennett and Sharpley in \cite[Theorem 5]{BS86} (cf. \cite[Lemma 7.5]{BS88}) and they used the idea of Lorentz and Shimogaki from \cite{LS71}, the so called ``pushing mass" technique. On the other hand, the original proof of Calder\'on's theorem \cite[Theorem 2.10, p. 114]{BS88} is based on the Hardy-Littlewood-P\'olya result. So that we will prove the monotone refinement of the Hardy-Littlewood-P\'olya theorem and use it in the alternative proof of Theorem \ref{CaldMitjMP}.

\proof[Proof of Theorem \ref{thmHLP}]
Proof of the first statement will be done by induction argument with respect to dimension $n$. The statement for $n=1$ is evident. Let $n>1$, assume the claim is true for all $k=1,...,n-1$ and let $0\leq a,b \in \mathbb{R}^n$ be both nonincreasing with $b\preceq a$. If $a_1=b_1$, then $a'=(a_2,a_3,...,a_n)$ and $b'=(b_2,b_3,...,b_n)$ satisfies $b'\preceq a'$. Thus by induction hypothesis there is $(n-1)\times (n-1)$ doubly stochastic and monotone matrix $B'$ such that $b'=B'a'$. Moreover, the matrix 
\[
B=\left[
\begin{array}{cc}
1&0
\\0&B'
\end{array}
\right]
\] 
is also doubly stochastic monotone and $Ba=b$. If $a_1>b_1$, we will find a doubly stochastic and monotone matrix $A'$ such that $b\preceq A'a$ but with $[A'a]_1=b_1$. Then it is enough to apply the previous step to vectors $b,A'a$ and the desired matrix will be $A=BA'$. Therefore we need only to find $A'$ like above. 

Suppose $a_1>b_1$. Because $\sum_{i=1}^na_i=\sum_{i=1}^nb_i$, we can find $k\leq n$ such that
\[
\frac{1}{k}\sum_{i=1}^ka_i\leq b_1 <\frac{1}{k-1}\sum_{i=1}^{k-1}a_i.
\] 
Consider the function 
\[
f(\eta)=\frac{\eta}{k-1}\sum_{i=1}^{k-1}a_i+(1-\eta)a_k
\]
for $0\leq \eta \leq 1$. 
We have $f(1)=\frac{1}{k-1}\sum_{i=1}^{k-1}a_i> b_1$ and $f(\frac{k-1}{k})=\frac{1}{k}\sum_{i=1}^{k-1}a_i+\frac{1}{k}a_k\leq b_1$. Therefore there is a solution $\gamma$ of the equation $f(\eta)=b_1$ that belongs to the interval $[\frac{k-1}{k},1)$.
Another words, the equation
\[
b_1=\sum_{i=1}^n\lambda_ia_i
\]
has a solution $\lambda_i=\frac{\gamma}{k-1}$ for $i=1,...,k-1$, $\lambda_k=(1-\gamma)$ and $\lambda_i=0$ for $i>k$. Moreover, the sequence $(\lambda_i)$ is nonincreasing because $\frac{k-1}{k}\leq\gamma\leq 1$. Then the matrix $A'$, that we are looking for, is 
\[
A'=\left[
\begin{array}{ccccccc}
\frac{\gamma}{k-1}&\dots&\frac{\gamma}{k-1}&1-\gamma&0&\dots&0
\\\vdots&\ddots&\vdots&\vdots&\vdots&&\vdots
\\\frac{\gamma}{k-1}&\dots&\frac{\gamma}{k-1}&1-\gamma&0&&0
\\1-\gamma&\dots&1-\gamma&\sigma &0&\dots&0
\\0&\dots&0&0&1&\dots&0
\\\vdots& &\vdots&\vdots&&\ddots&
\\0&\dots&0&0&0&\dots&1
\end{array}
\right],
\] 
where the entry $\sigma=1-(k-1)(1-\gamma)$ is on the position $(k,k)$. 
In fact, we have
\[
\sum_{i=1}^jb_i\leq jb_1=\sum_{i=1}^j[A'a]_i\ \ {\rm for}\ \ j\leq k-1,
\]
\[
\sum_{i=1}^jb_i\leq \sum_{i=1}^ja_i=\sum_{i=1}^j[A'a]_i\ \ {\rm for}\ \ j\geq k,
\]
so that $b\preceq A'a$. Moreover, $A'$ is evidently monotone and doubly stochastic. 

Note that above we have used strongly equality $\sum_{i=1}^na_i=\sum_{i=1}^nb_i$ therefore the second case, when we have only  $b\prec a$, has to be treated in another way.
We divide the interval $(0,n]$ of natural numbers into the collection of intervals $(0,i_1],(i_1,i_2],...,(i_{k-1},i_k=n]$ in such a way that $\delta _{j}a$ and $b$ restricted to a given interval $(i_{j-1},i_j]$ satisfies stronger relation $\preceq $ with some constant $\delta _{j}\leq 1$, i.e. $b\chi_{(i_{j-1},i_j]}\preceq \delta _{j}a\chi_{(i_{j-1},i_j]}$. Then it will be enough to find doubly stochastic and monotone matrices $A_j$ for each $(i_{j-1},i_j]$ and to define $A$ as
\[
A=\left[
\begin{array}{cccc}
\delta_1A_1& \dots &0
\\ \vdots &\ddots  & \vdots
\\ 0& \dots & \delta_kA_k
\end{array}
\right].
\]
It remains therefore to find the mentioned factorization of $(0,n]$ and a sequence $(\delta_j)$. 
There is nothing to do when $b\preceq a$. Suppose $b\prec a$. 
Put $i_0=0$ and define for $j=1,2,...$ 
$$\delta_{j+1}=\max_{k\leq n} \frac{\sum_{i=i_{j}+1}^kb_i}{\sum_{i=i_{j}+1}^ka_i}\ {\rm and }\  i_j=\max \{k\leq n:\sum_{i=i_{j}+1}^kb_i=\delta_{j+1} \sum_{i=i_{j}+1}^ka_i\}.
$$
We apply it until $i_{j}=n$ for some $j$. Denote such a last $j$ as $k$. Then  sequences $(\delta_{j})_{j=1}^k$ and $((i_{j-1},i_j]_{j=1}^k)$ are the desired ones. 
\endproof

\begin{remark}
To see that the above method is essentialy different than the Lorentz-Shimogaki ``pushing mass'' technique applied by Bennett and Sharpley in their proof (\cite[Lemma 3]{BS86}), let us consider vectors $f=(1,1,1,1)$ and $g=(2,1,1,0)$. Of course, $f\preceq g$. Applying the ``pushing mass" algorithm from \cite{BS86} we need three steps, i.e. $g_0=(3/2,3/2,1,0)$, $g_1=(4/3,4/3,4/3,0)$, $g_2=f=(4/3,4/3,4/3,0)$. On the other hand, our method produce the desired matrix immediately (we mean the inductive step from the first part of the above proof). It is also worth to mention that the maximal number of steps in our method is $n$ while the ``pushing mass'' technique may require more than $n$ steps (in general less than $2n$). In fact, consider  $f=(12/10,11/10,1,8/10)$ and $g=(2,1,1,0)$. Consequently, the steps of ``pushing mass'' method are $g_0=(3/2,3/2,1,0)$, $g_1=(4/3,4/3,4/3,0)$, $g_2=(12/10,12/10,12/10,4/10)$, $g_3=(12/10,11/10,11/10,6/10)$, $g_4=f$, and so $4$ steps were not enough to get $f$.  
\end{remark}
\begin{remark}
There is also another important connection of the proof of Theorem \ref{thmHLP}  with the classical results. Namely, careful reading ensures that it may be regarded as a discrete version of the ``cutting corners" method of Arazy and Cwikel (\cite[figure 1]{AC84}). It is however disappointing that their method cannot be applied to get the $(L^p,L^{\infty})$ version of Theorem \ref{CaldMitjMP} in case $p>1$ because then respective operators $S$ from \cite[pp. 258--260]{AC84} are not monotone. 
\end{remark}

Theorem \ref{CaldMitjMP} was stated without a proof as a corollary from its ``simple--function" version in \cite{BS86} and \cite{BS88}. It seems to be intuitively evident, but we explain a little more careful that the standard ``limit" argument preserves monotonicity of resulting operator. In order to do it we start with sketching main steps of the proof of  \cite[Proposition 2.9, p. 110]{BS88} and then we will complement the required explanation. 

\proof[Proof of Theorem \ref{CaldMitjMP}]
Let $0\leq f,g \in L^1+L^{\infty}$ be both nonincreasing with $g\prec f$. Suppose first that
\begin{equation}
g=\sum_{k=1}^n b_k\chi_{A_k},\label{simpleg}
\end{equation} 
where $A_k=[(k-1)d,kd)$ for some $d>0$. Define operators 
$G:L^1+L^{\infty}\rightarrow \mathbb{R}^{n}$ and $H:\mathbb{R}^{n}\rightarrow L^1+L^{\infty}$ as 
\begin{equation}
G : h \mapsto \Big(\frac{1}{\mu(A_k)}\int_{A_k}h\, d\mu \Big)_{k=1}^{n}\label{defG}
\end{equation}
and
\begin{equation}
H: (a_k)_{k=1}^{n} \mapsto \sum_{k=1}^{n} a_k \chi_{A_k}.\label{defH}
\end{equation}
Then the composition $HG$ is just an averaging operator and  $g=HG g$. Moreover, $g\prec f$ implies $G g\prec G f$ as well as $HG g\prec HGf$. We apply Theorem \ref{thmHLP} to find substochastic monotone matrix $B$ such that $G g= BG f$. Then put $T=HBG$ so that
$$
g=HG g= HBG f=Tf
$$ 
and $T$ is monotone because each of its components evidently is. Thus we have proved the thesis in the case when $g$ is of the form (\ref{simpleg}). Let now $g$ be arbitrary with $g\prec f$. 
We find a sequence $(g_m)$ such that $g_m\rightarrow g$ $\mu$ - a.e. and  $g_m\prec f$ with $g_m$ being nonincreasing for each $m$. Moreover, we may and do assumption that each $g_m$ is like in (\ref{simpleg}) (for some $n$ and $d$ depending on $m$).  Then one can apply the previous part to find sequence of substochastic monotone operators $T_m$ satisfying $T_m(f)=g_m$. Following the proof of \cite[Proposition 2.9, p. 110]{BS88}, it remains to define the set function 
$\nu_h:\Sigma |_E \rightarrow \mathbb{R}$, where $E$ is a measurable set with finite measure, by
\[
\nu_h(F)=\eta(\int_FT_mhd\mu)
\]
for $\eta$ being a fixed Banach limit. 
Then, exactly as in \cite{BS88}, we conclude that $\nu_h$ is absolutely continuous with respect to $\mu$ (restricted to $E$) and consequently $Th$ on $E$ is defined to be the Radon-Nikodym derivative of $\nu_h$ with respect to $\mu$. Since $E$ was arbitrary and by uniqueness of the Radon-Nikodym derivative  (up to a set of measure $0$) one can ''glue together" all parts to define $Th$ on the whole semiaxis. Also $Tf=g$ and the only that we need to explain more carefully is monotonicity of $T$. Let $0\leq h \in L^1+L^{\infty}$ be nonincreasing. Notice first that by the Lebesgue theorem $\frac{1}{\mu(F_n(t))}\int_{F_n(t)}Thd\mu\rightarrow Th(t)$ for almost all $t>0$, where $F_n(t)=(t-1/n,t+1/n)\cap \mathbb{R}_+$. Denote by $Z$ the set of such $t$. Choose $0<t_0<t_1$ from $Z$. Then for each $m$ and each $n$ such that $F_n(t_0)\cap F_n(t_1)=\emptyset$
\[
\frac{1}{\mu(F_n(t_0))}\int_{F_n(t_0)}T_mhd\mu\geq \frac{1}{\mu(F_n(t_1))}\int_{F_n(t_1)}T_mhd\mu
\] 
by monotonicity of $T_mh$. Finally, by positivity of functional $\eta$ we conclude that also  
\begin{equation*}
\frac{1}{\mu(F_n(t_0))}\int_{F_n(t_0)}Thd\mu\geq \frac{1}{\mu(F_n(t_1))}\int_{F_n(t_1)}Thd\mu 
\end{equation*}
for each such $n$ and thus $Th(t_0)\geq Th(t_1)$, which means that $Th$ is nonincreasing in the sense of (\ref{defmon}). 
\endproof

\begin{remark}
Dmitriev in \cite{Dm81} was considering the so called {\it positively $K$-monotone interpolation}, where having positive $f\in X_0+X_1,g\in Y_0+Y_1$ with $K(t,g;Y_0,Y_1)\leq K(t,f;X_0,X_1)$ for all $t>0$ one asks if there is a positive operator $T:(X_0,X_1)\rightarrow (Y_0,Y_1)$ satisfying $Tf=g$. Therefore Theorem \ref{CaldMitjMP} may be read as: the couple $(L^1,L^{\infty})$ satisfies ``monotone'' version of the above property, i.e. in place of positivity we require also monotonicity of $T$, assuming that $f,g$ are nondecreasing.   
\end{remark}

\section{Calder\'on couples}
We start with the following theorem which is the second of the main two steps toward the proof that $(\widetilde{L^1},L^{\infty})$ is a Calder\'on couple.
\begin{theorem} \label{thm:mapstotilde}
Let $f \in \widetilde{X}$ where $X$ is a Banach function space. Then for each $q>1$ there is a linear operator $T$ such that $Tf=\widetilde{f}$ and $\|T\|_{\widetilde{X}\rightarrow \widetilde{X}}\leq q$.
\end{theorem}
\proof
Let $f \in \widetilde{X}$, $f\not = 0$ and  $q>1$.
Since $h\mapsto \sign(f)h$ acts boundedly with norm one in each Banach function space, we may assume that $0\leq f$. 
For each $n\in \mathbb{Z}$ define
\[
A_n'=(q^{-n},q^{-n+1}].
\]
and
\[
A_n=(\tilde{f})^{-1}(A_n').
\]
Because $\tilde{f}$ is nonincreasing and right-continuous, such defined $A_n$ are either empty or are left-closed intervals. Moreover, we can choose a nondecreasing sequence of real numbers $(a_n)_{n=m_0}^{m_1}$ in such a way that 
each nonempty $A_n$ is of the form
\[
A_n=[a_{n},a_{n+1}),
\] 
or $A_n=[a_n,\infty)$ for some $n=m_1$. Let us say a few words about this sequence.   
We have three possibilities on the ``right side" of the sequence $(a_n)_{n=m_0}^{m_1}$. If $0<\lim_{t\rightarrow \infty}\tilde{f}(t)\in A'_k$ for some $k$, then we put $A_k=[a_k, \infty)$, $m_1=k$ and we additionally define $a_{m_1+1}=\infty$. The second case with finite ${m_1}$ occurs when there is $0<c<\infty$ such that  $\tilde{f}(t)=0$ for each $t\in [c,\infty)$ but $0<\lim_{t\rightarrow c^-}\tilde{f}(t)\in A'_k$ for some $k$. Then we understand that $m_1=k$, $A_k=[a_k,c)$ and $a_{m_1+1}=c$. In the remaining case, $m_1=\infty$ and for each $n_0$ there is $n_1>n_0$ such that $a_{n_1}>a_{n_0}$. The situation on the ''left" of $(a_n)$ is easier because either  $0<a_n\downarrow 0$ with $n\rightarrow -\infty$ and we put then $m_0=-\infty$, or $a_k=0$ for some $k$ and we put $m_0=k$ for the biggest such $k$ (notice that if  $f \not= 0$ then $a_k \not= 0$ for some $k$). Note also that $A_n$ may be nonempty only for $m_0\leq n \leq m_1$, where we understand $n<\infty$ if $m_1=\infty$ and $-\infty<n$ if $m_0=-\infty$. Of course, it may happen that $A_n$ is empty for some $n$ between $m_0$ and $m_1$, then we understand that $a_{n}=a_{n+1}$.   

We can now proceed with construction of the desired operator. 
For each $m_0\leq n \leq m_1$ limits 
$$
\tilde f(a_{n+1}^-)=\lim_{t\rightarrow a_{n+1}^-}\tilde{f}(t)
$$
are well defined and finite by monotonicity of $\tilde f$. Therefore, for each $m_0\leq n\leq m_1$ such that $A_n\not=\emptyset$ there exists a sequence of sets $(B^n_k)$ of finite, positive measure such that 
\[
B^n_k\subset \big(a_{n+1}-\frac{1}{k},\infty\big)
\] 
and
\begin{equation}
\tilde f\big(a_{n+1}-\frac{1}{k}\big) - \frac{1}{k}\leq f(t) \leq \tilde f\big(a_{n+1}-\frac{1}{k}\big) \label{indeks}
\end{equation}
for each $t\in B^n_k$ and $k\in \mathbb{N}$ satisfying $a_{n+1}-\frac{1}{k}\geq 0$, with the only exception in case $a_{m_1+1}=\infty$ (i.e. when $\tilde f(\infty)>0$), in which we just set 
\[
B^{m_1}_k\subset (ka_{m_1},\infty)
\] 
satisfying 
\begin{equation*}
\tilde f(\infty) - \frac{1}{k}\leq f(t) \leq \tilde f(\infty).
\end{equation*}
Let $\eta$ be a Banach limit. Then for each $n$ like above
\begin{equation}
\tilde f(a_{n+1}^-) =\lim_{k\rightarrow \infty}\frac{1}{\mu(B^n_k)}\int_{B^n_k}fd\mu =\eta_{k\uparrow}\Big( \frac{1}{\mu(B^n_k)}\int_{B^n_k}fd\mu \Big).
\end{equation}
We define the operator $S$ on $\widetilde X$ 
by the formula
\[
Sh=\sum_{n=m_0}^{m_1} \lambda_n (h)\chi_{A_n},
\]
where 
\[
 \lambda_n (h)= \eta_{k\uparrow} \Big( \frac{1}{\mu(B^n_k)}\int_{B^n_k}hd\mu\Big)
\]
for those $n$ with nonempty $A_n$ and $\lambda_n (h)=0$ for the remaining $n$ from the scale. 
Such $S$ is evidently linear and we will show that 
\begin{equation}\label{okreslonosc}
|Sh|\leq \tilde h.
\end{equation}
In fact, choose nonempty $A_n$ and let $t\in A_n$. We can find $l$ such that $B^n_k\subset [t,\infty)$ for all $k\geq l$. Thanks to it and by property (iii) of the Banach limit we get  
\begin{eqnarray*}
|Sh(t)|&=&| \lambda_n (h)|=|\eta_{k\uparrow}\big( \frac{1}{\mu(B^n_k)}\int_{B^n_k}hd\mu\big)| \\&=&|\eta_{k\uparrow}\big( \frac{1}{\mu(B^n_{k+l})}\int_{B^n_{k+l}}hd\mu\big)|\\
&\leq& \eta_{k\uparrow} \big(\frac{1}{\mu(B^n_{k+l})}\int_{B^n_{k+l}}|h|d\mu\big)\leq 
 {\rm ess}\sup_{s\geq t}|h(s)|=\tilde h(t).
\end{eqnarray*}
It means that $\|S\|_{\widetilde{X}\rightarrow\widetilde{X}}\leq 1$.
Therefore we have found the main part of the desired operator $T$. The second part will be simpler, just a multiplication operator $M_v:h\mapsto vh$ whose symbol $v$ is given by 
\[
v=\sum_{n=m_0}^{m_1} \frac{1}{\tilde f(a_{n+1}^-)} \tilde f \chi_{A_n}.
\]
It is well defined because $ \tilde f(a_{n+1}^-)>0$ for each $m_0\leq n\leq m_1$. Moreover, if $t\in A_n$ then
\[
\frac{\tilde f(t)}{\tilde f(a_{n+1}^-)}\leq \frac{q^{-n+1}}{q^{-n}}=q.
\]
Then $\|v\|_{L^{\infty}}\leq q$ which implies that  $\|M_v\|_{\widetilde{X}\rightarrow \widetilde{X}}\leq q$ (see for example \cite{MP89}). 
The proof is finished now, since 
\[
M_vSf=v(\sum_{n=m_0}^{m_1} \lambda_n (f)\chi_{A_n})=(\sum_{n=m_0}^{m_1} \frac{1}{\tilde f(a_{n+1}^-)} \tilde f \chi_{A_n})(\sum_{n=m_0}^{m_1} \tilde f(a_{n+1}^-)\chi_{A_n})=\tilde f.
\]
\endproof

The $K$-functional for the couple $(\widetilde{L^1},L^{\infty})$ was already calculated by Sinnamon in \cite{Si91} and it is given by the formula
\begin{equation}
K(t,f;\widetilde{L^1},L^{\infty})=K(t,\tilde f;L^1,L^{\infty})=\int_0^t \tilde f(s)ds. \label{functionalKK}
\end{equation}
It is also not difficult to see that the first equality may be extended to all couples $(X,L^{\infty})$ with $X$ being an arbitrary Banach function space.
\begin{proposition}\label{Kfunctionalfalka}
Let $X$ be a Banach function space and $f\in \widetilde{X}+L^{\infty}$. Then 
\[
K(t,f;\widetilde{X},L^{\infty})=K(t,\tilde f;X,L^{\infty}).
\]
\end{proposition}
\proof
Without loss of generality we may assume that $0\leq f\in \widetilde{X}+L^{\infty}$. It is enough to notice that 
\[
\widetilde{[(f-a)_+]}=(\tilde f-a)_+,
\]
where, as usually, $g_+=g\chi_{\{s:g(s)\geq 0\}}$.
Then 
\begin{eqnarray*}
K(t,f;\widetilde{X},L^{\infty})&=& \inf_{a>0}\{\|(f-a)_+\|_{\widetilde{X}}+at\}\\
&=& \inf_{a>0}\{\|(\tilde f-a)_+\|_{X}+at\}\\
&=& K(t,\tilde f;X,L^{\infty}).
\end{eqnarray*}
\endproof

\begin{remark}\label{Kfunctionalgeneral}
If we agree to replace equality $K(t,f;\widetilde{X},L^{\infty})=K(t,\tilde f;X,L^{\infty})$ by equivalence, then it holds in much more general setting. Namely, equivalence $K(\cdotp,f;\widetilde{X},\widetilde{Y})\approx K(\cdotp,\tilde f ;X,Y)$ for symmetric spaces $X,Y$ is a straightforward consequence of the general property that the operation $X\mapsto \widetilde{X}$ commutes with the Calder\'on-Lozanovskii construction, i.e. $\varphi (\widetilde{X},\widetilde{Y})= \widetilde{\varphi (X,Y)}$ with norms satisfying inequalities  
$$\|x\|_{\widetilde{\varphi (X,Y)}} \leq \|x\|_{\varphi (\widetilde{X},\widetilde{Y})}\leq C \|x\|_{\widetilde{\varphi (X,Y)}},$$
where $1\leq C\leq 2$ (see \cite{LM14c} and note that this equivalence holds also for some non symmetric spaces).  
In particular, taking $\varphi (u,v)=u+v$ and understanding that the space $tY$ contains the same elements as $Y$ with $\|f\|_{tY}=t\|f\|_{Y}$, we obtain
\[
\widetilde{X+tY}=\widetilde{X}+t\widetilde{Y},
\]
for each $t>0$. Consequently, 
\begin{eqnarray*}
K(t,\tilde g;X,Y)&=& \|\tilde g\|_{X+tY} = \|g\|_{\widetilde{X+tY}}  \leq \|g\|_{\widetilde{X}+t\widetilde{Y}}\\&\leq& 2\|g\|_{\widetilde{X+tY}}=2\|\tilde g\|_{X+tY}=2K(t,\tilde g;X,Y).
\end{eqnarray*}
and, since $K(t,g;\widetilde{X},\widetilde{Y})=\|g\|_{\widetilde{X}+t\widetilde{Y}}$, the claim follows. Notice that the above equivalences may be also deduced from \cite{CEP99} (cf. \cite{CC05}). 
\end{remark}

\begin{lemma}\label{MPOper}
If an operator $T:X\rightarrow Y$ is monotone, then $T:\widetilde{X}\rightarrow \widetilde{Y}$ with $\|T\|_{\widetilde{X}\rightarrow \widetilde{Y}}\leq \|T\|_{X\rightarrow Y}$.
\end{lemma}
\proof
Let $f\in \widetilde{X}$. Then, by monotonicity of $T$, we have
\[
\widetilde{T(f)}=\widetilde{|T(f)|}\leq \widetilde{T(|f|)}\leq \widetilde{T(\tilde f)}= T(\tilde f),
\]
which means that
\[
\|Tf\|_{\widetilde{Y}}= \|\widetilde{Tf}\|_{Y}\leq  \|T(\tilde f)\|_{Y}\leq \|T\|_{X\rightarrow Y}\|\tilde f\|_{X}=\|T\|_{X\rightarrow Y} \| f\|_{\widetilde{X}}.
\]
\endproof

We are now ready to state the main theorem of this section. 
\begin{theorem} \label{thm:monotCalderon}
The couple $(\widetilde{L^1},L^{\infty})$ is a Calder\'on couple.
\end{theorem}
\proof
Let $f,g \in \widetilde{L^1}+L^{\infty}$ with 
\[
K(t,g;\widetilde{L^1},L^{\infty})\leq K(t,f;\widetilde{L^1},L^{\infty})~{\rm for ~all}~ t>0.
\]
We will find $H$ satisfying $Hf=g$ according to the following scheme
\[
\begin{tikzcd}
f\arrow[mapsto]{r}{S} 
&\tilde{f}\arrow[mapsto]{r}{T}
&\tilde g\arrow[mapsto]{r}{W}
&g,
\end{tikzcd}
\]
where all of $S,T,W$ act boundedly from $(\widetilde{L^1},L^{\infty})$ into itself. 
Firstly, we find the last operator $W$ which is just multiplication by the function 
\[
\frac{g}{\tilde g}\leq 1,
\] 
where we understand $\frac{g(t)}{\tilde g(t)}=0$ when $\tilde g(t)=0$. 
In consequence, $\|W\|_{(\widetilde{L^1},L^{\infty})\rightarrow (\widetilde{L^1},L^{\infty})}=1$. 
Also operator $S$ is already known, because it is exactly the one from Theorem \ref{thm:mapstotilde}, let's say with the norm $\|S\|_{\widetilde{L^1+L^{\infty}}\rightarrow \widetilde{L^1+L^{\infty}}}=\gamma>1$. 
Notice that we apply Theorem \ref{thm:mapstotilde} for space $\widetilde{L^1}+L^{\infty}=\widetilde{L^1+L^{\infty}}$ (by  Proposition \ref{Kfunctionalfalka}) but then property (\ref{okreslonosc}) of the construction ensures also $\|S\|_{(\widetilde{L^1},L^{\infty})\rightarrow (\widetilde{L^1},L^{\infty})}\leq \gamma$.
It remains to find $T$. By Proposition \ref{Kfunctionalfalka}, the assumption  
\[
K(t,g;\widetilde{L^1},L^{\infty})\leq K(t,f;\widetilde{L^1},L^{\infty})
\]
means that 
\[
\int_0^t \tilde g(s)ds \leq \int_0^t \tilde f(s)ds,
\]
for all $t>0$. 
Therefore, applying Theorem \ref{CaldMitjMP} to functions $\tilde f,\tilde g$, we find a monotone operator $T$ such that $T\tilde f=\tilde g$ and $\|T\|_{(L^1,L^{\infty})\rightarrow (L^1,L^{\infty})}\leq 1$. Monotonicity of $T$ and Lemma \ref{MPOper} imply that also $\|T\|_{(\widetilde{L^1},L^{\infty})\rightarrow (\widetilde{L^1},L^{\infty})}\leq 1$. Finally, $H=WTS$ and the proof is finished with $\|H\|_{(\widetilde{L^1},L^{\infty})\rightarrow (\widetilde{L^1},L^{\infty})}\leq \gamma$. 
\endproof


Suppose we want to show that couples $(\widetilde{X_0},\widetilde{X_1}),(\widetilde{Y_0},\widetilde{Y_1})$ are relative Calder\'on couples. 
The proof of above theorem suggests the following point of view. Fix $g\in \widetilde{Y_0}+\widetilde{Y_1},f\in \widetilde{X_0}+\widetilde{X_1}$ with  $K(t,g;\widetilde{Y_0},\widetilde{Y_1})\leq K(t,f;\widetilde{X_0},\widetilde{X_1})$ for all $t>0$ and consider the scheme 
\[
\begin{tikzcd}
f\arrow[mapsto]{r}{A} 
&\tilde{f}\arrow[mapsto]{r}{S}
&\tilde g\arrow[mapsto]{r}{B}
&g\\
(\widetilde{X_0},\widetilde{X_1})\arrow{r}{A} 
&(\widetilde{X_0},\widetilde{X_1})\arrow{r}{S}
&(\widetilde{Y_0},\widetilde{Y_1})\arrow{r}{B}
&(\widetilde{Y_0},\widetilde{Y_1}).
\end{tikzcd}
\]   
Notice that once again existence of $A$ is a consequence of Theorem \ref{thm:mapstotilde} and $B$ is just a multiplication by function $g/\tilde g$. The only that misses is $S$. However, assumption on $f,g$ is, by Remark \ref{Kfunctionalgeneral}, equivalent (up to some constant) with $K(t,\tilde g;Y_0,Y_1)\leq K(t,\tilde f;X_0,X_1)$ for all $t>0$. Therefore, if we have proved that for a given positive, nonincreasing functions $h\in X_0+X_1,w\in Y_0+Y_1$ with $K(t,w;Y_0,Y_1)\leq K(t,h;X_0,X_1)$ for all $t>0$, there is a positive monotone operator $S:(X_0,X_1)\rightarrow (Y_0,Y_1)$ with $Sh=w$ then, thanks to Lemma \ref{MPOper}, we would have the desired operator $S:(\widetilde{X_0},\widetilde{X_1})\rightarrow  (\widetilde{Y_0},\widetilde{Y_1})$.  
According to this observation and using Dmitriev's results \cite{Dm74} we can generalize the main theorem to the Lorentz space setting. 

Recall that the Lorentz space $\Lambda_{\varphi}$ is defined by
$$
\Lambda_{\varphi}=\{f\in L^0: \| f \|_{\Lambda _{\varphi}}=\int f^{\ast}(t) d\varphi(t) < \infty \},
$$
where $\varphi$ is a concave, positive and increasing function on $[0,\infty)$ with  $\varphi(0^+)=0$ and $\varphi(\infty)=\infty$ (cf. \cite{BS86}, \cite{KPS82}). We get the following monotone version of Dimitriev's theorem \cite{Dm74}. 

\begin{theorem}\label{thmDimitriev}
Let $0\leq f \in \Lambda_{\varphi}+L^{\infty}$ and $0\leq g\in L^1+L^{\infty}$ be both nonincreasing and such that 
\[
K(t,g;{L^1},L^{\infty})\leq K(t,f;{ \Lambda_{\varphi}},L^{\infty}) ~ {\rm for ~ all} ~ t>0.
\]
Then there exists a monotone $S:(\Lambda_{\varphi},L^{\infty})\rightarrow (L^1,L^{\infty})$ such that $Sf=g$. 
\end{theorem}
\proof Similarly as in the proof of Theorem \ref{thm:monotCalderon} we sketch important steps of Dimitriev's proof to demonstrate where the monotone modification is necessary.
Let $0\leq f \in \Lambda_{\varphi}+L^{\infty}$ and $0\leq g\in L^1+L^{\infty}$ satisfy our assumptions. Suppose first that 
\begin{equation}
g=\sum_{k=1}^n b_k\chi{A_k},\label{simpleg2}
\end{equation} 
where $A_k=[(k-1)d,kd)$ for some $d>0$. For $h\in  \Lambda_{\varphi}+L^{\infty}$ define an operator $D$ in the following way
\[
Dh(z)= \sum_{k=1}^n \frac{1}{d}\int_{\varphi^{-1}((k-1)d)}^{\varphi^{-1}(kd)}h(t)d\varphi(t)\chi_{A_k}.
\]
Then $D$ is positive, $\|D\|_{(\Lambda_{\varphi},L^{\infty})\rightarrow (L^1,L^{\infty})}\leq 1$ and, especially, is monotone (cf. \cite[pages 529--530]{Dm74}). Moreover, $Df$ is nonincreasing and $g\prec Df$ so we can apply Theorem \ref{CaldMitjMP} to find monotone substochastic operator $T$ with $TDf=g$. Since $T$ and $D$ are monotone,  also $TD$ is monotone and we take just $S=TD$. Now, let a sequence $(g_n)$ consists of functions of the form (\ref{simpleg2}) and be such that $g_n\uparrow g$ a.e.. For each $g_n$ we find monotone operator $S_n$ like above. Then, once again following Dmitriev's explanation, we find the desired $S$ as an accumulation point of the sequence  $(S_n)$ with respect to the weak operator topology $\Gamma$ (see Appendix below), thanks to the result of Sedaev \cite{Se71}. Taking a subsequence if necessarily, we may assume that $(S_n)$ tends to $S$. In particular, it means that for each measurable $A\subset [0,\infty)$ with $\mu (A)<\infty$ and each function $h\in \Lambda_{\varphi}+L^{\infty}$ we have
\[
\int_AS_nhd\mu\rightarrow \int_AShd\mu {\rm ~as~} n\rightarrow \infty.
\]
Therefore, one can explain monotonicity of $S$ like in the proof of Theorem \ref{CaldMitjMP}. 
\endproof

\begin{corollary}
The couple $(\widetilde{\Lambda_{\varphi}},L^{\infty})$ is a relative Calder\'on couple with respect to $(\widetilde{L^1},L^{\infty})$.
\end{corollary}

Note that it is not necessary to follow the way described above to conclude that some couple of the form $(\widetilde{X},\widetilde{Y})$ is a Calder\'on couple. In fact, a straightforward application of Theorem 14 from \cite{AC12} to Theorem \ref{thm:monotCalderon} gives such a result for $(\widetilde{L^p},L^{\infty})$, although we know nothing about monotone operators in this case. In fact, we see that $\widetilde{(|f|^p)}=(\tilde f)^p$, which means that $p$-convexification $(\widetilde{X})^p$ of $\widetilde{X}$ is exactly $\widetilde{(X^p)}$. Recall that $p$-convexification $X^p$ ($p\geq 1$) of a Banach function space $X$ is $X^p=\{f:|f|^p\in X\}$ with the norm $\|f\|_{X^p}= \||f|^p\|^{\frac{1}{p}}_{X}$.

\begin{theorem}
For $1\leq p <\infty$, $(\widetilde{L^p},L^{\infty})$ is a Calder\'on couple.
\end{theorem}

\begin{remark}
All the above results, except 
Remark \ref{Kfunctionalgeneral}, remain true when the underlying measure space $(\mathbb{R}_+,\mu)$ will be replaced by $([0,1],\mu)$ with the Lebesgue measure $\mu$. 
\end{remark}

According to the above considerations the following question seems to be of interest. 

\begin{problem}
Let $X,Y$ be symmetric spaces such that $(X,Y)$ is a Calder\'on couple. Let $0\leq f,g \in X+Y$ be both nonincreasing and such that 
\[
K(t,g;X,Y)\leq K(t,f;X,Y) ~ {\rm for ~ all} ~ t>0.
\]
Does there exist monotone operator $T$ acting on the couple $(X,Y)$ with $Tf=g$, or is $(\widetilde{X},\widetilde{Y})$ is a Calder\'on couple? As we have already seen,  $(\widetilde{L^p},L^{\infty})$ is a Calder\'on couple but we do not know if there is a monotone operator like above. Notice that respective operators from Lorentz-Shimogaki, Cwikel and Arazy-Cwikel papers are not monotone. Also the proof of Theorem 14 from \cite{AC12} says nothing about this, because it is  based on the lattice version of Hahn-Banach extension theorem. 
\end{problem}

{\bf Acknowledgements} The author is very grateful to Professor Lech Maligranda for valuable remarks and advices which allowed to improve the paper.

\section{Appendix}

The paper of Sedaev \cite{Se71} is not easy to acquire and is in Russian while, one the other hand, Dmitriev's explanation is quite abbreviated. Because of these facts, just for the sake of convinience, we recall Sedaev's result and explain how it is applied in Theorem \ref{thmDimitriev}.  

We introduce some special notion after Sedaev, while the remaining terminology is the standard one like in books \cite{BL76}, \cite{BS88} or \cite{KPS82}. 
Let $\overline{X}=(X_0,X_1),\overline{Y}=(Y_0,Y_1)$ be two couples of compatible  Banach spaces. Further, let 
\[
\Gamma \subset (Y_0+Y_1)^*,
\] 
\[
G=\{T:\|T\|_{(X_0,X_1)\rightarrow (Y_0,Y_1)}\leq 1 \},
\]
\[
U_i=\{y\in Y_i:\|y\|_{Y_i}\leq 1\}, ~ i=1,2.
\]
Then $\sigma(Y_0+Y_1,\Gamma)$ means the weak topology on $Y_0+Y_1$ restricted to $\Gamma$. Similarly, $\Gamma$-topology on $L(\overline{X},\overline{Y})$ is the weak operator topology restricted to $\Gamma$. Denote also, after Sedaev, 
\[
a=\inf_{x\in X_0\cap X_1}\frac{\|x\|_{X_1}}{\|x\|_{X_0}},b=\sup_{x\in X_0\cap X_1}\frac{\|x\|_{X_1}}{\|x\|_{X_0}}.
\]
\begin{theorem}[Sedaev 1971]
$G$ is $\Gamma$-compact in $L(\overline{X},\overline{Y})$ if and only if
\begin{itemize}
\item[(i)] $U_0\cap cU_1$ is $\sigma(Y_0+Y_1,\Gamma)$-closed in $Y_0+Y_1$ for each $a\leq c\leq b$,
\item[(ii)] if $X_0\cap X_1$ is not dense in $X_0$ ($X_1$) then $U_0$  ($U_1$) is $\sigma(Y_0+Y_1,\Gamma)$-closed in $Y_0+Y_1$,
\item[(iii)] there is a couple $(X,Y)$ such that 
\begin{itemize}
\item[(a)] $X$ is an interpolation space for $\overline{X}$ and $E$ is dense in $X_0+X_1$,
\item[(b)] $Y$ is an intermediate space for $\overline{Y}$ and the unit ball $U$ of $Y$ is  $\sigma(Y_0+Y_1,\Gamma)$-compact.
\end{itemize}
\end{itemize}
\end{theorem}

The family of operators $(S_n)$ from the proof of Theorem \ref{thmDimitriev} clearly belong to $G$, where $(X_0+X_1)=(\Lambda_{\varphi},L^{\infty}), (Y_0,Y_1)=(L^1,L^{\infty})$. Moreover, we choose 
\[
\Gamma = \{f\in L^{\infty}:\mu(\supp(f))<\infty  \}.
\] 
Of course, $\Gamma \subset L^1\cap L^{\infty}=(L^1+L^{\infty})'$.
We need only to explain that such couples and $\Gamma$ satisfies assumptions of the above theorem. To prove $(i)$ and $(ii)$ we will show that both $U_0$ and $U_1$ are $\sigma(L^1+L^{\infty},\Gamma)$-closed in $L^1+L^{\infty}$. Let $f\in L^1+L^{\infty}$ be such that $f\not\in U_0$. This means that 
\[
\int_{\mathbb{R}_+} |f|d\mu>1+3\delta 
\] 
for some $\delta>0$.
Then there is $A\subset \mathbb{R}_+$, $\mu(A)<\infty$ such that 
\[
\int_A |f|d\mu>1+2\delta. 
\]
Set $g=\sign(f)\chi_A\in \Gamma$. Then $$\langle g,f\rangle=\int_{\mathbb{R}_+}  gfd\mu=\int_A |f|d\mu >1+2\delta$$ and $$V=\{h\in  L^1+L^{\infty}: |\langle g,f-h\rangle|<\delta\}$$ is an $\sigma(L^1+L^{\infty},\Gamma)$-open neighbourhood of $f$. Moreover, since $g\in U_1$, it follows that $|\langle g,h\rangle|\leq 1$ for each $h\in U_0$ and consequently $V\cap U_0=\emptyset $, which means that $U_0$ is $\sigma(L^1+L^{\infty},\Gamma)$-closed in $L^1+L^{\infty}$. Consider now $U_1$. Analogously as before, let $f\in L^1+L^{\infty}$ be such that $f\not\in U_1$. This means that there is $A\subset [0, \infty)$, $\mu(A)<\infty$ such that $|f|\chi_A>(1+\delta)\chi_A$. This time put $g=\frac{\sign(f)}{\mu(A)}\chi_A\in \Gamma$. Then $$\langle g,f\rangle=\int_{\mathbb{R}_+}  gfd\mu>1+\delta$$ and we can proceed as before, because $g\in U_0$.   

To see that the point $(iii)$ is satisfied we set $X=Y=L^{\infty}$. Point $(a)$ is satisfied because simple functions are dense in $\Lambda_{\varphi}$. It remains to notice that $U_1$ is $\sigma(L^1+L^{\infty},\Gamma)$-compact. But $U_1$ is $\sigma(L^{\infty},L^1)$-compact in $L^{\infty}$ thanks to Alaoglu theorem. It is then also compact in the weaker topology $\sigma(L^{\infty},\Gamma)$ (since $\Gamma\subset L^1$) and so also in $L^1+L^{\infty}$ with topology $\sigma(L^1+L^{\infty},\Gamma)$.    


\end{document}